# Energy Conserving Schemes for the Simulation of Musical Instrument Contact Dynamics[†]

Vasileios Chatziioannou[*]

*Institute of Music Acoustics (IWK), University of Music and performing Arts Vienna, Anton-von-Webern-Platz 1, 1030 Vienna, Austria*

Maarten van Walstijn

*Sonic Arts Research Centre, Queen's University Belfast, BT7 1NN Belfast, Northern Ireland*

**Abstract**

Collisions are an innate part of the function of many musical instruments. Due to the nonlinear nature of contact forces, special care has to be taken in the construction of numerical schemes for simulation and sound synthesis. Finite difference schemes and other time-stepping algorithms used for musical instrument modelling purposes are normally arrived at by discretising a Newtonian description of the system. However because impact forces are non-analytic functions of the phase space variables, algorithm stability can rarely be established this way. This paper presents a systematic approach to deriving energy conserving schemes for frictionless impact modelling. The proposed numerical formulations follow from discretising Hamilton's equations of motion, generally leading to an implicit system of nonlinear equations that can be solved with Newton's method. The approach is first outlined for point mass collisions and then extended to distributed settings, such as vibrating strings and beams colliding with rigid obstacles. Stability and other relevant properties of the proposed approach are discussed and further demonstrated with simulation examples. The methodology is exemplified through a case study on tanpura string vibration, with the results confirming the main findings of previous studies on the role of the bridge in sound generation with this type of string instrument.

## 1. Introduction

When studying the vibrational behaviour of musical instruments or other sounding objects, collisions are often encountered. These can occur either in a confined space [1] (e.g. hammer-string interaction, mallet impacts) or in a more distributed manner [2], such as the coupling between the snares and the membrane of a snare drum. The former can usually be modelled in lumped form, suppressing the computation of the interaction forces to a single point, whereas the latter require considering variations along spatial coordinates. In both cases, the impactive interaction represents an important nonlinear element that is closely linked to the characteristics and/or the expressive control of the instrument [3].

[*]Corresponding author Tel.: +43 1 71155 4313.
*Email addresses:* chatziioannou@mdw.ac.at (Vasileios Chatziioannou), m.vanwalstijn@qub.ac.uk (Maarten van Walstijn)

Contact modelling has undergone extensive study, including various proposed time-stepping methods for the simulation of vibro-impact phenomena. A fundamental distinction between two different classes herein stems from whether or not an interpenetration is allowed between the contacting objects. A perfectly rigid contact involves a non-penetration condition of the form

$$y_a - y_b \geq 0, \tag{1}$$

where a moving body with position $y_a$ collides with another body located below it, at $y_b$. The use of Lagrange multipliers is common in conjunction with this approach and finite element simulations [4]. When penetration (which can be equivalently considered as the compression of the impacting objects) is allowed, a repelling force can be defined using a penetration function

$$P(\chi) = h(\chi)\chi \tag{2}$$

where $\chi = y_b - y_a$ and $h(\chi)$ denotes the Heaviside step function. This is referred to as a penalty approach [5], the validity of which is subject to constraints on impact velocity and penetration level [6]. Using either of the above methodologies, existence and uniqueness of solutions has been proven only for a small number of special cases [7].

The penalty approach commonly appears in musical acoustics problems in the form of a one-sided power law [8, 9], where, starting from Hertz's contact law, the impact force takes the form

$$f(\chi) = k_c \lfloor \chi^\alpha \rfloor \tag{3}$$

where $\lfloor \chi^\alpha \rfloor = h(\chi)\chi^\alpha$, $k_c$ is a stiffness coefficient and the power law exponent $\alpha \geq 1$ depends on the local shape of the contact surface. These parameters are often determined empirically, and good agreement with measurements has been found for several cases [1, 10–12]. A version of Eq. (3) with impact friction is possible in the form of the Hunt-Crossley model [13–15], which has also found use in various other engineering fields (e.g. robotics [16]).

When simulation is required to solve collision problems, the power law needs to be incorporated into a numerical formulation. Most of the relevant time-stepping methods found in the musical acoustics and sound computing literature are based on finite differences [8, 10] or closely related methods such as the trapezoidal rule [14], the Newmark-beta method [11], or Verlet integration [17, 18]. In distributed settings, discretisation is sometimes performed after first casting the linear part of the problem in modal form [19, 20]. While many successful simulation results have been obtained, and stability can even be shown for some specific cases or under specific assumptions (see e.g. [21]), the formulation of a more general class of numerical schemes for impact modelling is still considered as an open and difficult problem [8, 17]. Sound-related collisions have also been simulated with digital waveguides [22–25], wave digital filters [26] and hybridisations thereof [27, 28]. Stability in such wave-variable models is generally analysed and controlled through passivity of the individual scattering units. However, provably stable formulations of this type for distributed impact governed by Eq. (3) are yet to appear.

Seeking a more rigorous numerical treatment of vibro-impact problems, the mathematical physics literature suggests that two distinct directions can be taken. One approach is to design a method such that the total energy is maintained, leading to energy conserving schemes [29, 30]. Alternatively, one may choose to preserve another invariant of the physical system, the symplectic structure [31], thus deriving symplectic numerical schemes, some of which have recently been applied to musical instrument sound synthesis [17, 18]. It has been shown that in general only one of the above properties can be preserved [32]. Although both approaches can establish the stability of a numerical algorithm, symplectic schemes allow an oscillating energy that can distort the amplitude of lossless systems. Such schemes are therefore particularly suited to the study of families of trajectories and long-term behaviour of dynamical systems, while the use of energy preserving schemes has been indicated as more suitable for oscillatory problems [33].

Energy based methods for constructing time-stepping algorithms can be generally classified into two categories, as explained in [34]; those yielding schemes that attempt to conserve an energy-like, positive definite quantity and those that aim to conserve the actual energy of the system at each time step. In strict physical terms, only the latter can be specified as '*energy conserving schemes*'. For the sake of clarity, the former energy based methods will be referred to as '*energy methods*' in the remaining of this manuscript. This distinction is made more clear in Section 3.2, where conserved numerical quantities for both types of methods are compared.

The past decade has seen a substantial uptake of such '*energy methods*' in application to nonlinear problems encountered in musical acoustics, most notably by Bilbao [8]. The merits of this approach—which requires identifying a numerical counterpart of the system Hamiltonian—come to the fore whenever the studied vibrational behaviour is intrinsically nonlinear, as is the case for various percussion instruments [35, 36]. However, one-sided power laws such as that in Eq. (3) are non-analytic functions of the phase-space variables, making it far from straightforward to derive schemes for which an invariant, numerical energy-like quantity exists [8]. The present authors propose to address this by first reformulating the system in its Hamiltonian form [30, 37], and discretise this rather than Newton's equations of motion, in order to construct an '*energy conserving scheme*'. This approach bears resemblance to that taken by Greenspan in discretising nonlinear lumped systems [38], and a similar strategy has recently been applied in extended form by Chabassier et al. to the simulation of nonlinear string vibrations [39]. Preliminary results of the application to collision problems have been reported by the current authors in [40], focusing on simulation of a point mass interacting with a barrier. The present study extends this to distributed interaction by modelling of a stiff string colliding with a (nearly) rigid obstacle, which has direct application to simulation and sound synthesis of various string instruments.

Since the stability analysis is carried out by virtue of the preservation of an invariant energy, frictional forces are initially neglected in the derivation and analysis of the numerical schemes. A way to introduce damping in compatibility with all presented algorithms is shown towards the end of this paper, which is organised as follows. Section 2 outlines the proposed methodology through the case of a simple one-mass system involving collisions, with specific focus on proving convergence and energy conservation, further supported by numerical examples. A similar treatment is then adopted in Section 3 for distributed collisions, which involves discussion of additional aspects such as matrix formulation, dispersion, and boundary conditions. Section 4 presents a case study on the interaction of a tanpura string with a curved bridge and Section 5 evaluates the main findings within the context of musical acoustics and sound synthesis.

## 2. Lumped contact

An elemental, frictionless model is defined, in which a gravity-sensitive mass $m$ attached to a spring of constant $k$ is colliding with a rigid barrier positioned at $y = y_c$. Assuming the impact force to be of the form given in the previous section (with $\chi = y_c - y$), the motion of the mass is governed by

$$m\frac{d^2 y}{dt^2} = k_c \lfloor (y_c - y)^\alpha \rfloor - ky + mg_0 \tag{4}$$

where $g_0$ is the gravitational acceleration (taken negative). Considering an energetic description, the Lagrangian of the system governed by Eq. (4), defined as the difference between kinetic energy $T$ and potential energy $V$, is

$$L(y, \dot{y}) = \frac{1}{2}m\dot{y}^2 - \frac{k}{2}y^2 - \frac{k_c}{\alpha + 1}\lfloor (y_c - y)^{\alpha+1} \rfloor + mg_0 y. \tag{5}$$

Defining the conjugate momentum $p = \partial L/\partial \dot{y}$ and taking the Legendre transformation of the Lagrangian yields the Hamiltonian of the system

$$H(y, p) = \frac{p^2}{2m} + \frac{k}{2}y^2 + \frac{k_c}{\alpha + 1}\lfloor (y_c - y)^{\alpha+1} \rfloor - mg_0 y = T(p) + V(y). \tag{6}$$

This equals the total energy of the system and is constant in this case due to the absence of frictional or external (non-conservative) forces. The corresponding Hamilton's equations of motion are [37]

$$\frac{dy}{dt} = \frac{\partial H(y, p)}{\partial p} = \frac{\partial T(p)}{\partial p} \tag{7a}$$

$$\frac{dp}{dt} = -\frac{\partial H(y, p)}{\partial y} = -\frac{\partial V(y)}{\partial y}. \tag{7b}$$

## 2.1. Numerical formulation

Amongst the various possible ways of solving the system numerically, a notable example regularly employed in musical and speech acoustics is the trapezoidal rule, which yields an implicit scheme [14]. Another possibility is to factorise the collision force term in the right-hand side of Eq. (4) into $k_c(y_c - y) \cdot \lfloor (y_c - y)^{\alpha-1} \rfloor$ and apply an averaging operator to the first term of this product, while approximating the left-hand side term of Eq. (4) with a centred difference term. This yields a scheme that allows an explicit update form [8]. Both schemes are unconditionally stable in the absence of the nonlinear collision term. However in neither is the energy of the numerical system preserved through transitions of $y$ across $y_c$, and for the explicit scheme the stability within compressed-state simulation phases is ensured only for specific integer exponent values [8].

A more general treatment follows from discretising the Hamiltonian rather than the Newtonian description. If $y^n$ denotes the value of variable $y$ at time $n\Delta t$, with $\Delta t$ being the sampling interval, employing mid-point derivative approximations for all terms in Eq. (7) yields

$$\frac{y^{n+1} - y^n}{\Delta t} = \frac{T(p^{n+1}) - T(p^n)}{p^{n+1} - p^n} \tag{8a}$$

$$\frac{p^{n+1} - p^n}{\Delta t} = -\frac{V(y^{n+1}) - V(y^n)}{y^{n+1} - y^n}. \tag{8b}$$

Setting $q^n = p^n \Delta t/(2m)$ and $\xi = \Delta t^2/(2m)$ allows writing scheme (8) as

$$y^{n+1} - y^n = q^{n+1} + q^n \tag{9a}$$

$$q^{n+1} - q^n = -\xi \frac{V(y^{n+1}) - V(y^n)}{y^{n+1} - y^n}. \tag{9b}$$

Solving Eq. (9) is facilitated by defining the auxiliary variable

$$s = y^{n+1} - y^n = q^{n+1} + q^n \tag{10}$$

which gives

$$q^{n+1} = s - q^n, \qquad y^{n+1} = y^n + s. \tag{11}$$

Substituting into Eq. (9b) gives a nonlinear function in $s$

$$F(s) = \xi \frac{V(y^n + s) - V(y^n)}{s} + s - 2q^n = 0 \tag{12}$$

with

$$\lim_{s \to 0} F(s) = \xi V'(y^n) - 2q^n \tag{13}$$

where $V'$ signifies taking the derivative of $V$ with respect to displacement, hence there is no singularity in $F(s)$. In comparison, the application of the trapezoidal rule results in an implicit scheme of exactly the same form, but with the nonlinear equation to be solved defined as

$$F_{\text{tr}}(s) = \xi \frac{V'(y^n + s) + V'(y^n)}{2} + s - 2q^n = 0 \tag{14}$$

whereas application of the implicit midpoint rule, which is a symplectic method [31], yields

$$F_{\text{mr}}(s) = \xi V'\left(\frac{2y^n + s}{2}\right) + s - 2q^n = 0. \tag{15}$$

Note that equations (12), (14) and (15) are equivalent for any potential $V$ that defines a linear force $f = \frac{\partial V}{\partial y}$, for example when $k_c = 0$. However for any nonlinear force function, such as that defined for contact, these expressions—and therefore the resulting numerical schemes—are distinct.

## 2.2. Numerical solution

Solution of the numerical scheme relies on finding a physically correct root of Eq. (12), which is subsequently used to update $y$ and $q$ at each time step using Eq. (11). For robustness of the scheme, it is useful to know that the function (12) has a unique solution, which can be shown as follows. From the definition of $F(s)$ it follows that

$$\frac{dF}{ds} = 1 + \xi \frac{s\,V'(y^n + s) - V(y^n + s) + V(y^n)}{s^2} \tag{16}$$

with $\lim_{s \to 0} \frac{dF}{ds} = 1 + \frac{\xi}{2} V''(y^n)$. Demonstrating that $\frac{dF}{ds} \geq 1$ proves that $F(s)$ always has a single root. This is equivalent to showing that

$$V(y^n + s) \leq V(y^n) + s\,V'(y^n + s) \tag{17}$$

which holds by definition $\forall\, y_n \in \mathbb{R}$, since the potential $V$ is a convex function of $y$. Hence a unique solution of Eq. (12) can be found using the Newton-Raphson method, which is globally convergent for a convex function [41]. This can be shown to hold for $F$ by writing

$$\frac{d^2F}{ds^2} = \xi \frac{Q(s)}{s^3}, \tag{18}$$

where

$$Q(s) = s^2 V''(y^n + s) - 2s V'(y^n + s) + 2V(y^n + s) - 2V(y^n). \tag{19}$$

Since $\frac{dQ}{ds} = s^2 V'''(y^n + s) \geq 0$ for any potential with convex first derivative (which holds $\forall \alpha \geq 1$), $Q(s)$ is a monotonically increasing function going through the origin, which from Eq. (18) implies that $\frac{d^2F}{ds^2} \geq 0$. This also holds in the limit, in that $\lim_{s \to 0} \frac{d^2F}{ds^2} = \frac{\xi}{3} V'''(s) \geq 0$. The number of iterations required for the solution of Eq. (12) can be kept low (typically below 6) by using the previous value of $s$ as the initial guess.

## 2.3. Conservation of energy

The principal advantage of the presented scheme is that it inherently conserves the total energy of the numerical system, which is readily demonstrated by rewriting Eq. (8) as

$$\frac{1}{\Delta t}(y^{n+1} - y^n)(p^{n+1} - p^n) = T(p^{n+1}) - T(p^n) \tag{20a}$$

$$\frac{1}{\Delta t}(y^{n+1} - y^n)(p^{n+1} - p^n) = -V(y^{n+1}) + V(y^n) \tag{20b}$$

and substituting by parts, which yields

$$T(p^{n+1}) + V(y^{n+1}) = T(p^n) + V(y^n) \quad \Rightarrow \quad H(y^{n+1}, p^{n+1}) = H(y^n, p^n). \tag{21}$$

This states that energy is conserved across each time step. Note that neither the above conservation proof nor the proof of global convergence to a unique solution relies on any assumptions regarding the range of the parameter values.

## 2.4. Accuracy

Beyond energy analysis, the immediate next question to explore is how well the scheme approximates the continuous-time model. While standard finite difference procedures [42] may be used to show that the scheme is of second order accuracy, further insight can be obtained by inspecting the manner in which the inherent approximation errors manifest themselves in the simulation results. For the case of linear oscillation ($k_c = 0$, $k > 0$), the system has a natural frequency of oscillation ($\omega_0 = k/m$), which is however not preserved in the numerical model. This can be established by substituting a single-frequency test solution of the form $y^n = e^{s_a \Delta t}$, $q^n = B e^{s_a \Delta t}$, into Eq. (9a). Using that $y^{n-1} = z^{-1} y^n$, where $z = e^{j\omega_d \Delta t}$, it follows that any continuous-domain complex frequency $s_a = j\omega_a$ maps to the discrete-domain complex frequency $z$ as:

$$s_a = B\frac{1 - z^{-1}}{1 + z^{-1}} = jB\tan(\omega_d \Delta t/2) \tag{22}$$

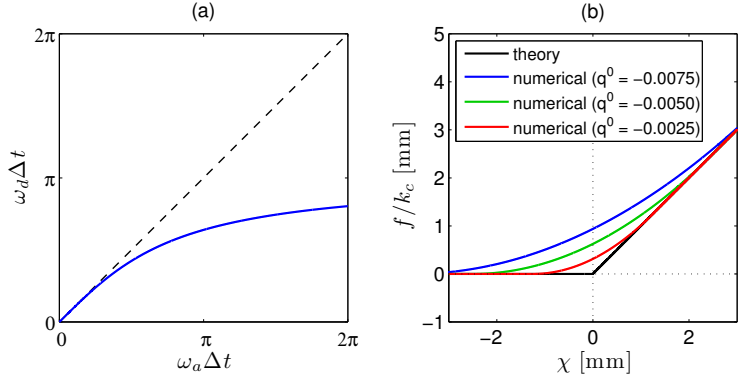

Figure 1: Inspection of the numerical approximation error. (a) Mapping between continuous-domain frequency $\omega_a$ and discrete-domain frequency $\omega_d$ for linear oscillation [$k_c = 0$]. (b) Effective repelling force for differently sized impact momenta, with $\alpha = 1$.

where $B = 2/T$. The above relationship is associated with the use of the bilinear transform in digital filter design [43]; indeed, as can be seen from expanding Eq. (12) and (14), for $k_c = 0$ the scheme in Eq. (9) is equivalent to that obtained with the trapezoidal rule. The associated warping of the frequency axis is shown in Fig. 1(a). The resonance shift may be pre-compensated for by setting $k$ to a higher value such that the numerical model possesses the correct resonance frequency, but such pre-warping has less significance in a nonlinear setting, and in addition does not extend readily to distributed problems. The approximation error can nevertheless generally be reduced by decreasing the time step $\Delta t$, at the cost of more computational effort.

The above frequency-domain analysis applies only to the linear dynamics of the system. For $k_c > 0$, the error committed within the nonlinear part of the system may be directly observed through the approximation of the repelling force. From Eq. (8b), the effective repelling force at time $t = (n + \frac{1}{2})\Delta t$ is

$$f = \frac{k_c}{\alpha + 1} \left[ \frac{\lfloor (y_c - y^n - s)^{\alpha+1} \rfloor - \lfloor (y_c - y^n)^{\alpha+1} \rfloor}{s} \right]. \tag{23}$$

Note that $s = q^{n+1} + q^n$ can be thought of as twice the mid-point value $q^{n+1/2}$, thus representing a normalised measure of impact momentum. Fig. 1(b) shows a zoom-in of the (normalised) effective repelling force ($f/k_c$) for three different impact momentum values, as directly evaluated from Eq. (23) against the mid-point compression $\chi^{n+1/2} = y_c - (y^{n+1} + y^n)/2$, and comparing to the corresponding theoretical term $\lfloor (\chi^{n+1/2})^\alpha \rfloor$. As can be seen, the scheme effectively smoothes the curve around $\chi = 0$. For $\alpha = 1$, this leads to a continuously differentiable force function, whereas the original force function was not differentiable at $\chi = 0$. For arbitrary $\alpha \geq 1$, the numerical model effectively replaces a force function of class $C^{\alpha-1}$ with one of $C^\alpha$. The discrepancy between the effective repelling force and its theoretical counterpart decreases rapidly with decreasing impact momentum size (by $1/s$), and the scheme converges to Eq. (3) in the limit:

$$\lim_{s \to 0} f = k_c \lfloor \chi^\alpha \rfloor. \tag{24}$$

Given that $s \to 0$ when $\Delta t \to 0$, this also demonstrates that the numerical model is consistent with theory.

## 2.5. Numerical examples of lumped system simulations

Adopting the expression of $V$ from Eq. (6), the general form of $F(s)$ can be expanded to

$$F(s) = \left(1 + \xi \frac{k}{2}\right) s + 2\left(\xi \frac{k}{2} y^n - q^n\right) - \xi m g_0 + \frac{\beta_c}{\alpha + 1} \frac{\lfloor (y_c - y^n - s)^{\alpha+1} \rfloor - \lfloor (y_c - y^n)^{\alpha+1} \rfloor}{s} \tag{25}$$

where $\beta_c = \xi k_c$. The code used for generating the results presented here directly solves this equation.

The most basic model described by Eq. (4) involves a single collision of a point mass with a rigid barrier (in the absence of gravity) and can be formulated by setting $k = 0$ and $g_0 = 0$. This example places the focus on the nonlinear part of $F(s)$. Existence and uniqueness of solutions, as well as conservation of energy are inherited from the general model for this and any subsequent example. However, due to quantisation in finite-precision arithmetic, the Hamiltonian can, at best, be preserved to machine precision in implementations on digital processors. It is therefore of interest to observe the resulting energy error, expressed in terms of the deviation of $H^n = H(y^n, p^n)$ from the initial energy $H^0$, which in normalised form reads

$$e^n = \frac{H^n - H^0}{H^0}. \tag{26}$$

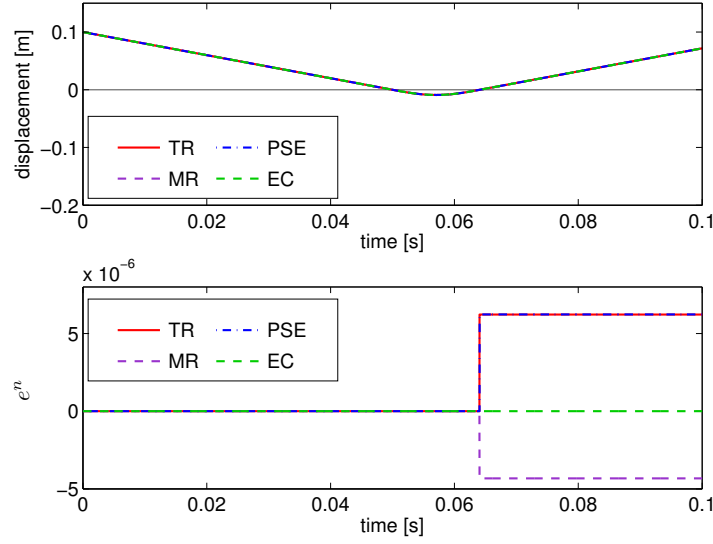

Figure 2: Collision of a mass ($m = 0.1$ kg) with a rigid barrier simulated using the presented energy conserving scheme (EC), the partially stable explicit scheme (PSE), the trapezoidal rule (TR) and the midpoint rule (MR). The stiffness is chosen as $k_c = 5000$ N m$^{-1}$ with $\alpha = 1$ and $f_s = \Delta t^{-1} = 44.1$ kHz. Top: mass displacement with initial position $y_0 = 0.1$ m and momentum $p_0 = -0.2$ kg m s$^{-1}$. Bottom: the energy error $e^n$.

It is worth noting here that quantisation generally results into a random-like signal $e^n$ which, if zero mean, will not cause an energy shift over time.

Fig. 2 compares the proposed energy conserving (EC) scheme (9) to the trapezoidal rule (TR), the midpoint rule (MR) and the partially stable explicit finite difference scheme (PSE) presented in [8], in terms of the simulated mass trajectory and the associated energy error size. The lower plot in Fig. 2 indicates that the PSE, TR and MR schemes can introduce energy jumps, which are observed here at the point of decoupling with the barrier. This artefact is avoided with the EC scheme, with the energy error barely exceeding machine precision levels.

The reduced scheme corresponding to a single mass-wall collision is described using only two parameters, namely $\alpha$ and $\beta_c$. In order to get a more complete view of the energy preservation properties of the proposed scheme, its performance is analysed across a range of values for these parameters, corresponding to different levels of interaction between the mass and the barrier. To ensure a meaningful comparison, the calculations are made independent of the collision duration and the initial energy of the system, using the following energy preservation metric:

$$\mathcal{P} = \sum_{n=n_1}^{n_2} \frac{|H^{n+1} - H^n|}{(n_2 - n_1 + 1)H^0} \tag{27}$$

where the collision occurs in the interval $[n_1, n_2]$. $\mathcal{P}$ can be thought of as the mean energy deviation per sample during the contact period, thus excluding periods during which energy deviations are expected to be negligible. As depicted in Fig. 3 the preservedness is only mildly dependent on the model parameters, and structurally retains very low values.

Let us now consider a ball falling under gravity and bouncing on the floor (at $y_c = 0$). This can be modeled by setting $g_0 = -9.81$ m s$^{-2}$ and $k = 0$. Fig. 4 (left) shows the results of such a simulation for $\alpha = 3.5$. In accordance with the conservation of energy in the absence of losses, the ball keeps bouncing back to its initial height.

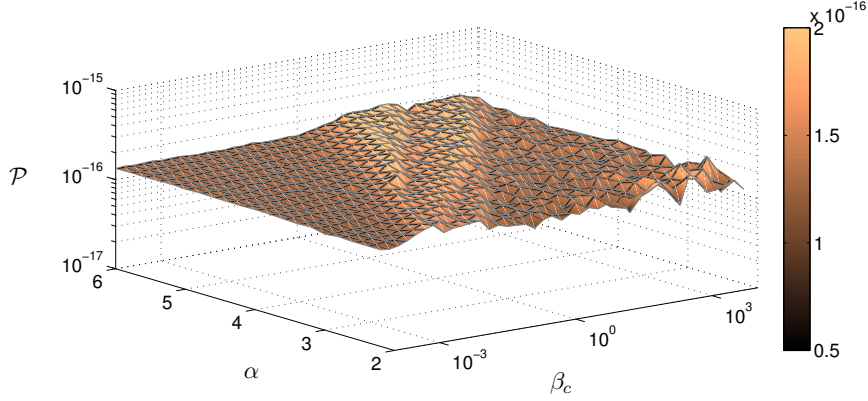

Figure 3: Simulation results of the energy preservation metric (27) as a function of $\alpha$ and $\beta_c$.

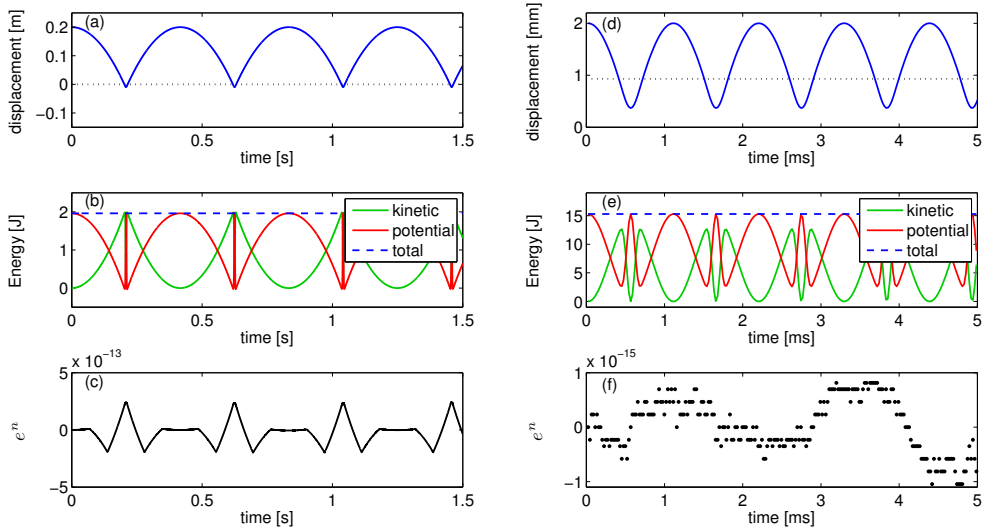

Figure 4: Left: Simulation of a lossless bouncing ball under a gravitational force with $k_c = 10^{11}$ and $\alpha = 3.5$: (a) Displacement (b) the corresponding energy components and (c) the energy error $e^n$. Right: Simulation of a lossless oscillating unit mass attached to a spring of stiffness $k = (2\pi 440)^2$ N m$^{-1}$. A repelling force becomes active when $y < 0.93$ mm, following a quadratic power law with $k_c = 2.5 \times 10^{10}$: (d) Displacement (e) the corresponding energy components and (f) the energy error $e^n$.

The system becomes more relevant to musical acoustics when the moving element can store potential energy, facilitating oscillatory behaviour. In the lumped model this is effected by setting the stiffness $k$ to a positive value. For an initial displacement value greater than the barrier position $y_c$, the scheme will now simulate a series of periodic interactions. Fig. 4 (right) shows an example for $\alpha = 2$. As can be seen, the repetitive impacts do not cause an accumulative energy shift, and the energy is conserved to machine precision. This was observed for a large number of simulations with different parameters and long simulation times.

## 2.6. Aliasing

Due to the heavy nonlinear character of the contact force as a function of displacement, any oscillation of the mass with barrier contact involves the generation of a series of overtones, which can potentially lead to significant aliasing in the simulation. To exemplify this, a series of simulations were performed, without gravity, for exponentially increasing $k$ values. This means that the 'system resonance' (i.e. the frequency of oscillation of the discrete model in the absence of collision) and the aforementioned overtones increase with each increase in $k$. Fig. 5(a) shows an example of the magnitude spectrum of the resulting momentum signal as a function of frequency for $\Delta t = 1/44.1$ms, plotted against the system resonance frequency. The aliasing terms can clearly be identified as the mirrored frequency components. This problem can be addressed by oversampling; Fig. 5(b) shows the spectrum for the same set of simulations but using $\Delta t = 1/176.4$ms, in which case the amplitudes of the overtones nearer the Nyquist frequency are small enough to avoid significant aliasing within the audible frequency range.

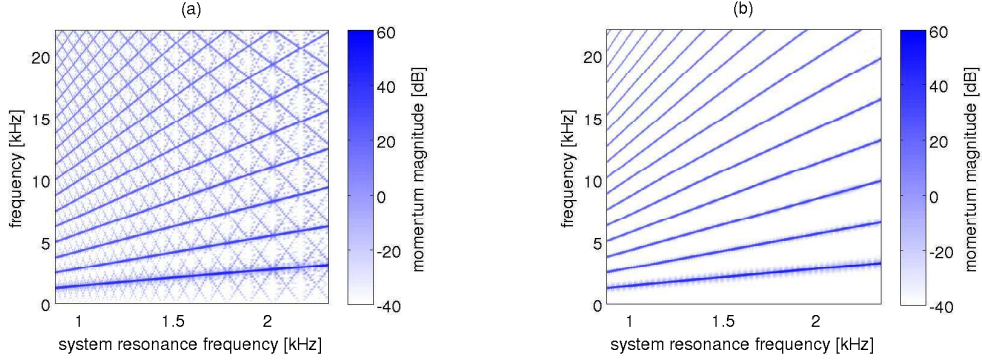

Figure 5: Magnitude spectrum of the calculated momentum signal as a function of frequency and system resonance with the sampling frequency set to (a) 44.1 kHz and (b) 176.4 kHz. The used parameters are $k_c = 2 \cdot 10^{10}$, $\alpha = 2.3$, $m = 0.001$ kg, $y_c = -0.05$m. The simulation was run for $k$ successively equalling $30000 \cdot 1.01^\iota$, $\iota = 1, 2 \ldots 200$, and the initial condition $y = 0.1$ m, $p = -0.1$ m/s at $n = 0$ was applied in each simulation. Each spectrum is computed by applying an FFT to 1 second of simulation output.

## 3. Distributed contact

The methodology to derive numerical schemes presented above can also be applied to distributed systems. Given the relevance of impactive interaction to string instrument vibrations, a useful representative case is that of a stiff string, the free transverse vibrations of which are governed by [44]

$$\rho A \frac{\partial^2 y}{\partial t^2} = \tau \frac{\partial^2 y}{\partial x^2} - EI \frac{\partial^4 y}{\partial x^4} \tag{28}$$

where $\rho$, $\tau$, and $A$ respectively are the mass density, tension, and cross-sectional area, while $E$ and $I$ denote the Young's modulus and the moment of inertia. This Euler-Bernoulli model can also represent a flexible string or ideal beam, by nulling $EI$ or $\tau$, respectively. Collisions with a distributed barrier can be included by adding a force density term of the form of Eq. (3) to the right hand side of Eq. (28), where $y_c(x)$ now represents the barrier profile. Hamilton's equations of motion for this system are then given by

$$\frac{\partial p}{\partial t} = \frac{\partial}{\partial x}\left(\frac{\partial \mathcal{H}}{\partial u}\right) - \frac{\partial^2}{\partial x^2}\left(\frac{\partial \mathcal{H}}{\partial v}\right) - \frac{\partial \mathcal{H}}{\partial y} \tag{29a}$$

$$\frac{\partial y}{\partial t} = \frac{\partial \mathcal{H}}{\partial p} \tag{29b}$$

where

$$\mathcal{H}(y, p, u, v) = \frac{1}{2}\frac{p^2}{\rho A} + \frac{1}{2}\tau u^2 + \frac{1}{2}EI v^2 + \frac{k_c}{\alpha + 1}\lfloor (y_c - y)^{\alpha+1} \rfloor = \mathcal{T}(p) + \mathcal{V}_\tau(u) + \mathcal{V}_s(v) + \mathcal{V}_c(y) \tag{30}$$

with

$$u = \frac{\partial y}{\partial x} \quad \text{and} \quad v = \frac{\partial^2 y}{\partial x^2}. \tag{31}$$

A derivation of these equations from a variational principle can be found in [45], where the simulation of an ideal, simply supported string is considered. $\mathcal{H}$ is the Hamiltonian density as a function of the local displacement $y$ and momentum density $p = \rho A\, \partial y/\partial t$ and the spatial derivatives $u$ and $v$. The total energy of the system with a string of length $L$ is the integral

$$H = \int_{x=0}^{L} \mathcal{H}(y, p, u, v)\, dx. \tag{32}$$

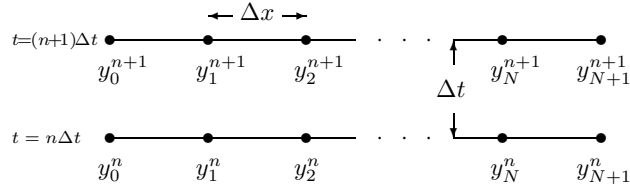

Figure 6: The discretised string at times $n\Delta t$ and $(n+1)\Delta t$.

*3.1. Numerical formulation*

Following the methodology of the previous section, mid-point derivative approximations are again employed to derive a numerical scheme, where $y_m^n$ now denotes the value of variable $y$ at position $x = m\Delta x$ and time $t = n\Delta t$, $\Delta x$ being the spatial sampling interval (see Fig. 6). Rewriting Eq. (29) using the separation of the energies in Eq. (30), Hamilton's equations are approximated by

$$\frac{p_m^{n+1} - p_m^n}{\Delta t} = \frac{\left\{\frac{\partial \mathcal{V}_\tau}{\partial u}\right\}_{m+\frac{1}{2}}^{n+\frac{1}{2}} - \left\{\frac{\partial \mathcal{V}_\tau}{\partial u}\right\}_{m-\frac{1}{2}}^{n+\frac{1}{2}}}{\Delta x} - \frac{\left\{\frac{\partial \mathcal{V}_s}{\partial v}\right\}_{m+1}^{n+\frac{1}{2}} - 2\left\{\frac{\partial \mathcal{V}_s}{\partial v}\right\}_{m}^{n+\frac{1}{2}} + \left\{\frac{\partial \mathcal{V}_s}{\partial v}\right\}_{m-1}^{n+\frac{1}{2}}}{\Delta x^2} - \frac{\mathcal{V}_c(y_m^{n+1}) - \mathcal{V}_c(y_m^n)}{y_m^{n+1} - y_m^n} \tag{33a}$$

$$\frac{y_m^{n+1} - y_m^n}{\Delta t} = \frac{\mathcal{T}(p_m^{n+1}) - \mathcal{T}(p_m^n)}{p_m^{n+1} - p_m^n} \tag{33b}$$

where, for example,

$$\left\{\frac{\partial \mathcal{V}_\tau}{\partial u}\right\}_{m+\frac{1}{2}}^{n+\frac{1}{2}} = \frac{\mathcal{V}_\tau(u_{m+\frac{1}{2}}^{n+1}) - \mathcal{V}_\tau(u_{m+\frac{1}{2}}^n)}{u_{m+\frac{1}{2}}^{n+1} - u_{m+\frac{1}{2}}^n} \tag{34}$$

approximates the partial derivative of $\mathcal{V}_\tau$ with respect to $u$. It is useful here to introduce the forward and backward space shift operators, through their action on $y_m^n$, as

$$\delta_+ y_m^n = \frac{y_{m+1}^n - y_m^n}{\Delta x}, \qquad \delta_- y_m^n = \frac{y_m^n - y_{m-1}^n}{\Delta x}. \tag{35}$$

Using the following approximations for $u$ and $v$

$$u_{m+\frac{1}{2}}^n = \delta_+ y_m^n, \qquad u_{m-\frac{1}{2}}^n = \delta_- y_m^n, \qquad v_m^n = \delta_+\delta_- y_m^n = \delta_\Delta y_m^n \tag{36}$$

a scheme centred at time $t = (n+1/2)\Delta t$ and position $x = m\Delta x$ is obtained

$$\frac{p_m^{n+1} - p_m^n}{\Delta t} = \frac{\tau}{2}\delta_\Delta(y_m^{n+1} + y_m^n) - \frac{EI}{2}\delta_\Delta^2(y_m^{n+1} + y_m^n) - \frac{k_c}{\alpha+1}\frac{\lfloor (y_{c_m} - y_m^{n+1})^{\alpha+1} \rfloor - \lfloor (y_{c_m} - y_m^n)^{\alpha+1} \rfloor}{y_m^{n+1} - y_m^n} \tag{37a}$$

$$\frac{y_m^{n+1} - y_m^n}{\Delta t} = \frac{1}{\rho A}\frac{p_m^{n+1} + p_m^n}{2} \tag{37b}$$

where $y_{c_m}$ denotes the location of the collision boundary at the position $x = m\Delta x$. In matrix form this can be written as

$$\mathbf{p}^{n+1} - \mathbf{p}^n = \phi \mathbf{D}_2\left(\mathbf{y}^{n+1} + \mathbf{y}^n\right) - \psi \mathbf{D}_4\left(\mathbf{y}^{n+1} + \mathbf{y}^n\right) - \varpi \mathbf{S}^{-1}\left(\lfloor(\mathbf{y}_c - \mathbf{y}^{n+1})^{\alpha+1}\rfloor - \lfloor(\mathbf{y}_c - \mathbf{y}^n)^{\alpha+1}\rfloor\right) \tag{38a}$$

$$\theta\left(\mathbf{y}^{n+1} - \mathbf{y}^n\right) = \mathbf{p}^{n+1} + \mathbf{p}^n. \tag{38b}$$

where $\mathbf{S} = \text{diag}(\mathbf{y}^{n+1} - \mathbf{y}^n)$ is a diagonal matrix,

$$\phi = \frac{\tau \Delta t}{2\Delta x^2}, \qquad \psi = \frac{EI \Delta t}{2\Delta x^4}, \qquad \varpi = \frac{k_c \Delta t}{\alpha + 1}, \qquad \theta = \frac{2\rho A}{\Delta t} \qquad (39)$$

and $\mathbf{y}^n$, $\mathbf{y}_c^n$ and $\mathbf{p}^n$ are column vectors holding displacement, barrier profile and momentum values. Under the assumption of simply supported boundary conditions on both ends of the system, these vectors hold the values of $N$ interior nodes on the string (i.e. from $y_1$ to $y_N$), and $\mathbf{D}_2$ then is an $N \times N$ tridiagonal matrix:

$$\mathbf{D}_2 = \begin{bmatrix} -2 & 1 & & 0 \\ 1 & \ddots & \ddots & \\ & \ddots & \ddots & 1 \\ 0 & & 1 & -2 \end{bmatrix} \qquad (40)$$

which implements the second spatial derivative of the string state, with $\mathbf{D}_4 = \mathbf{D}_2 \mathbf{D}_2$. Analogous to the lumped model derivations, it is convenient to rewrite the scheme using a scaled momentum variable $\mathbf{q}^n = \mathbf{p}^n/\theta$. Also substituting $\mathbf{D} = (\beta_4 \mathbf{D}_4 - \beta_2 \mathbf{D}_2)$ then gives

$$\mathbf{q}^{n+1} - \mathbf{q}^n = -\mathbf{D}\left(\mathbf{y}^{n+1} + \mathbf{y}^n\right) - \zeta \mathbf{S}^{-1}\left(\lfloor (\mathbf{y}_c - \mathbf{y}^{n+1})^{\alpha+1}\rfloor - \lfloor (\mathbf{y}_c - \mathbf{y}^n)^{\alpha+1}\rfloor\right) \qquad (41a)$$

$$\mathbf{y}^{n+1} - \mathbf{y}^n = \mathbf{q}^{n+1} + \mathbf{q}^n \qquad (41b)$$

where $\beta_2 = \phi/\theta$, $\beta_4 = \psi/\theta$ and $\zeta = \varpi/\theta$. Now setting

$$\mathbf{s} = \mathbf{y}^{n+1} - \mathbf{y}^n = \mathbf{q}^{n+1} + \mathbf{q}^n \qquad (42)$$

yields the nonlinear system of equations

$$\mathbf{F} = (\mathbf{I} + \mathbf{D})\mathbf{s} + 2(\mathbf{D}\mathbf{y}^n - \mathbf{q}^n) + \zeta \mathbf{S}^{-1}\left(\lfloor (\mathbf{y}_c - \mathbf{y}^n - \mathbf{s})^{\alpha+1}\rfloor - \lfloor (\mathbf{y}_c - \mathbf{y}^n)^{\alpha+1}\rfloor\right) = \mathbf{0} \qquad (43)$$

which is the distributed equivalent of Eq. (25) for $g_0 = 0$.

*3.2. Conservation of energy*

The total energy of the system can be calculated by integrating the energy densities along the length of the string, i.e.

$$H^n = \sum_{m=0}^{N+1} (\mathcal{T}_m^n + \mathcal{V}_m^n) \Delta x \qquad (44)$$

where $\mathcal{T}_m^n = \mathcal{T}(p_m^n)$ and

$$\mathcal{V}_m^n = \frac{\mathcal{V}_\tau(\delta_- y_m^n) + \mathcal{V}_\tau(\delta_+ y_m^n)}{2} + \mathcal{V}_s(\delta_\Delta y_m^n) + \mathcal{V}_c(y_m^n) \qquad (45)$$

represent kinetic and potential energy densities, respectively. This calculation involves so-called 'ghost nodes' lying just outside the spatial domain occupied by the string, which are eliminated by applying appropriate numerical boundary conditions (see Section 3.6). In matrix form the total energy is

$$H = b[\mathbf{q}^t \mathbf{q} + \mathbf{y}^t \mathbf{D} \mathbf{y} + \zeta \mathbf{1}^t \lfloor (\mathbf{y}_c - \mathbf{y})^{\alpha+1}\rfloor] \qquad (46)$$

with $\mathbf{1} = (1, \ldots, 1)^t$ and $b = 2\rho A \Delta x / \Delta t^2$. For simply supported ends this can be written in the more intuitive form

$$H^n = \Delta x \left[\frac{(\mathbf{p}^n)^t \mathbf{p}^n}{2\rho A} + \frac{\tau}{2\Delta x^2}(\mathbf{D}_1 \mathbf{y}^n)^t (\mathbf{D}_1 \mathbf{y}^n) + \frac{EI}{2\Delta x^4}(\mathbf{D}_2 \mathbf{y}^n)^t (\mathbf{D}_2 \mathbf{y}^n) + \frac{k_c}{\alpha + 1} \mathbf{1}^t \lfloor (\mathbf{y}_c - \mathbf{y}^n)^{\alpha+1}\rfloor\right] \qquad (47)$$

where

$$\mathbf{D}_1 = \begin{bmatrix} 1 & & & & & 0 \\ -1 & 1 & & & & \\ & -1 & 1 & & & \\ & & \ldots & \ldots & & \\ & & & & -1 & 1 \\ 0 & & & & & -1 \end{bmatrix} \quad (48)$$

is an $(N + 1) \times N$ matrix which implements taking the gradient of the string state, with $\mathbf{D}_1^t \mathbf{D}_1 = -\mathbf{D}_2$. Seeking to demonstrate the conservation of energy, we multiply the left hand side of Eq. (41a) with $(\mathbf{q}^{n+1} + \mathbf{q}^n)^t$ and the right hand side with $(\mathbf{y}^{n+1} - \mathbf{y}^n)^t$, which are equivalent terms by Eq. (41b). This yields

$$(\mathbf{q}^{n+1} - \mathbf{q}^n)^t(\mathbf{q}^{n+1} + \mathbf{q}^n) = -(\mathbf{y}^{n+1} - \mathbf{y}^n)^t \mathbf{D}(\mathbf{y}^{n+1} + \mathbf{y}^n) - \zeta(\mathbf{y}^{n+1} - \mathbf{y}^n)^t \mathbf{S}^{-1}\left(\lfloor (\mathbf{y_c} - \mathbf{y}^{n+1})^{\alpha+1}\rfloor - \lfloor (\mathbf{y_c} - \mathbf{y}^n)^{\alpha+1}\rfloor\right) \quad (49)$$

which, given that $\mathbf{D}$ is symmetric, can be written as

$$(\mathbf{q}^{n+1})^t\mathbf{q}^{n+1} + (\mathbf{y}^{n+1})^t\mathbf{D}\mathbf{y}^{n+1} + \zeta\mathbf{1}^t\left(\lfloor(\mathbf{y_c} - \mathbf{y}^{n+1})\rfloor^{\alpha+1}\right) = (\mathbf{q}^n)^t\mathbf{q}^n + (\mathbf{y}^n)^t\mathbf{D}\mathbf{y}^n + \zeta\mathbf{1}^t\left(\lfloor(\mathbf{y_c} - \mathbf{y}^n)\rfloor^{\alpha+1}\right). \quad (50)$$

Now multiplying by $b$ and using the definition of the numerical energy in Eq. (46), it follows that

$$H^{n+1} = H^n. \quad (51)$$

As explained in the Introduction, instead of an 'energy conserving scheme', different numerical schemes can be constructed using an 'energy method' that conserves an energy-like quantity, as proposed for example in [46]. In that case the conserved quantity, defined at time $t = (n + 1/2)\Delta t$ becomes

$$\widetilde{H}^{n+1/2} = (\mathbf{q}^{n+1/2})^t\mathbf{q}^{n+1/2} + (\mathbf{y}^{n+1})^t\mathbf{D}\mathbf{y}^n + \zeta\mathbf{1}^t\left(\frac{\lfloor(\mathbf{y_c} - \mathbf{y}^{n+1})\rfloor^{\alpha+1} + \lfloor(\mathbf{y_c} - \mathbf{y}^n)\rfloor^{\alpha+1}}{2}\right) \quad (52)$$

whereas the actual discrete system energy at time $t = (n + 1/2)\Delta t$ is equal to

$$H^{n+1/2} = (\mathbf{q}^{n+1/2})^t\mathbf{q}^{n+1/2} + (\mathbf{y}^{n+1/2})^t\mathbf{D}\mathbf{y}^{n+1/2} + \zeta\mathbf{1}^t\left(\lfloor(\mathbf{y_c} - \mathbf{y}^{n+1/2})\rfloor^{\alpha+1}\right) \quad (53)$$

and has a different potential energy term. This results in an oscillating system energy $H$, whereas an energy-like quantity $\widetilde{H}$ is conserved. Hence while allowing the definition of stable numerical schemes, such 'energy methods' do not intrinsically inherit the energy conservation property of the underlying model equations, neither do they replicate it; instead a conserved energy-like quantity has to be found for each specific case. Also note that, unlike $H^{n+1/2}$, $\widetilde{H}^{n+1/2}$ is not bound to be positive. This imposes a stability condition on the associated numerical schemes (see [46]).

### 3.3. Numerical solution

Eq. (43) can be solved for $\mathbf{s}$ using the multidimensional Newton method, which requires forming the Jacobian of $\mathbf{F}$

$$\mathbf{J} = \mathbf{I} + \mathbf{D} + \mathbf{C} \quad (54)$$

where $\mathbf{C}$ is a diagonal matrix with elements

$$\{c_{i,i}\} = \frac{\Delta t}{\theta}\frac{s_i \mathcal{V}'_c(y_i^n + s_i) - \mathcal{V}_c(y_i^n + s_i) + \mathcal{V}_c(y_i^n)}{s_i^2} \quad (55)$$

which in accordance with Eq. (17) is positive definite. From the energy expression (46) it follows that $\mathbf{D}$ and therefore also $\mathbf{J}$ are positive definite, which ensures the uniqueness of a root of Eq. (43) [41]. Singularities in both $\mathbf{F}$ and its Jacobian can be handled as in the lumped case. The update equation is $\mathbf{s} = \mathbf{s} - \mathbf{J}^{-1}\mathbf{F}$, where instead of forming the inverse matrix it is possible—and considerably more efficient—to solve a (band) linear system. Global convergence of the Newton method (for an arbitrary initial guess) is guaranteed for the componentwise convex function $\mathbf{F}$ when the Jacobian is an $M$-matrix [47], which holds for $\beta_4 = 0$. For non-zero stiffness the method is only locally convergent and a good initial guess is required, which in practice is always available through the previous value of $\mathbf{s}$. As such, convergence is typically achieved in fewer than 20 iteration steps.

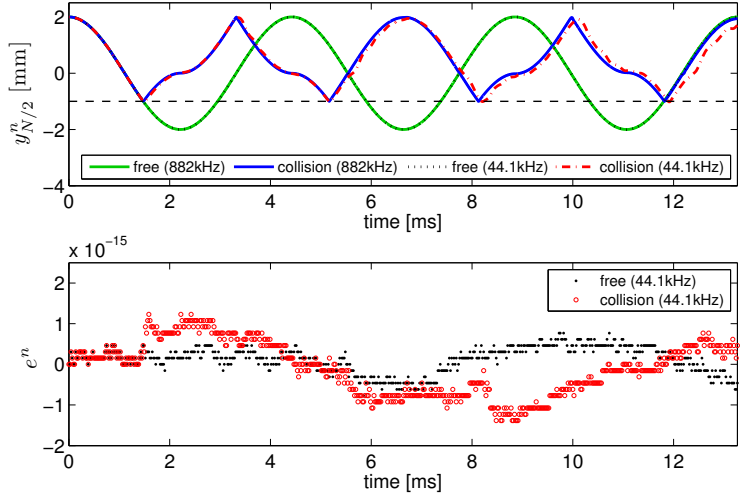

Figure 7: Simulation of an ideal string being free to vibrate or bouncing on a rigid obstacle, for the initial condition $y(x, 0) = 0.002 \sin(\pi x/L)$. The linear mass density of the string is set to $\rho A = 0.001$ kg m$^{-1}$ and the tension to $\tau = 100$ N, with $\Delta x = 0.007$ m and $f_s = \Delta t^{-1} = 44.1$ kHz. Top: mid-point string displacement. Bottom: numerical energy error.

### 3.4. Test simulations

To verify the correct behaviour of the distributed model, reference is made to an analytical result that compares the frequency of a free, flexible vibrating string with that of an impeded one [48]. This states that if a straight obstacle is placed halfway across the amplitude of the string vibration, the fundamental frequency of the free vibrating string will be 1.5 times the frequency of the impeded string [49]. This result is reproduced by a numerical simulation using scheme (41) for a 0.7 m long string (see Fig. 7), setting $k_c = 10^7$ in order to simulate a rigid obstacle with $\alpha = 1$. For comparison, the simulation was repeated with 20 times oversampling and using $k_c = 10^9$, which results in a very close approximation to the theoretical frequency ratio.

Fig. 8(a) shows the vibration of a stiff string bouncing on a curved surface located at one of its boundaries, setting $EI = 0.012$ N m$^2$. The nonlinear behaviour can be observed in the irregular exchange between the kinetic and the potential energy in Fig 8(b). Nonetheless the total energy remains constant, with inter-sample energy steps only occasionally exceeding machine precision levels (see Fig. 8(d)). Fig. 8(c) reveals that the interaction between the string and the boundary involves multiple impacts during contact periods, which result in the generation of high-frequency vibrations that are characteristic of string instruments with flat bridges. For example, the buzzing sound of a sitar is understood to stem from such multiple impacts [50]. Supplementary animations for both the above simulations, as well as the cantilever beam simulation of Section 3.6 are available at http://www.socasites.qub.ac.uk/mvanwalstijn/jsv14/.

### 3.5. Numerical dispersion

In the absence of the nonlinear collision term, the scheme in Eq. (41) reduces to that obtained by applying the Crank-Nicolson method [12, 42], which can be considered as the distributed version of the trapezoidal rule. Hence the frequency-domain approach of Section 2.4 applies again for $k_c = 0$, meaning that the mode frequencies of the string are warped according to the mapping in Fig. 1(a). At first sight, this appears to be a significant downside, as there are other finite difference schemes available that—for the linear case—introduce considerably less numerical dispersion (see, e.g. [8]). However, the effect is often only just audible at standard audio rates (e.g. 44.1kHz). More importantly, for the impactive case ($k_c > 0$) the collisions represent strongly nonlinear events driving the generation of high frequency components, which invariably results in aliasing effects (see also Section 2.6). Hence in practice the time step already has to be chosen about 2 to 4 times smaller than 1/44.1 ms to avoid detrimental effects on the simulation results.

The unconditional stability of the scheme represents a beneficial feature, in that the spatial and temporal stepsize can be chosen independently. Care must be taken however not to stretch the ratio between them too far. That is, the number of string modes modelled equals the number of string segments, thus increasing $N$ while holding $\Delta t$ constant amounts to compressing more and more modes into the finite simulation bandwidth according to the mapping in Fig. 1(a). It make sense therefore to not exceed a certain ratio to avoid rendering a string with very densely spaced modes at the high end of the spectrum. For example, for an ideal string with $c = \sqrt{\tau/\rho A}$ it is practical to not let the value of $\beta_2 = \frac{1}{4}(c\Delta t/\Delta x)^2$ exceed much beyond $1/4$, where setting $\beta_2 = 1/4$ corresponds to the numerical model possessing exactly all the modes with theoretical frequencies below Nyquist frequency.

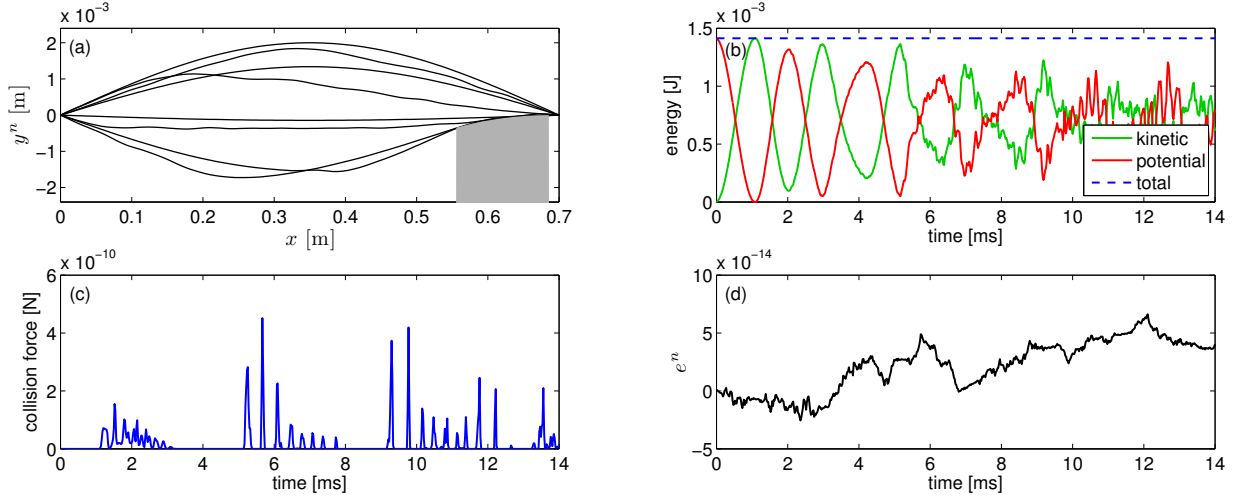

Figure 8: (a) Snapshots of a stiff string bouncing on a curved obstacle, (b) the energy components (c) the collision force due to string-obstacle interaction and (d) the energy error $e^n$.

Table 1: Numerical boundary conditions.

|  | clamped | simply supported | free |
|---|---|---|---|
| **left end** | $y_0^n = 0$ | $y_0^n = 0$ | $y_0^n = 0$ |
|  | $\delta_- y_0^n = 0$ | $\delta_\Delta y_0^n = 0$ | $\delta_- \delta_\Delta y_0^n = 0$ |
| **right end** | $y_{N+1}^n = 0$ | $y_{N+1}^n = 0$ | $\delta_\Delta y_{N+1}^n = 0$ |
|  | $\delta_+ y_{N+1}^n = 0$ | $\delta_\Delta y_{N+1}^n = 0$ | $\delta_+ \delta_\Delta y_{N+1}^n = 0$ |

### 3.6. Boundary conditions and damping

As seen in Eq. (43) and (46), the system and its numerical energy can be expressed in terms of the spatial matrix $\mathbf{D} = \beta_2 \mathbf{D}_2 - \beta_4 \mathbf{D}_4$, which for simply supported ends is entirely defined through the matrix $\mathbf{D}_2$ as given by Eq. (40). Other boundary conditions, such as free or clamped ends, can be effected simply by altering of and/or adding elements to $\mathbf{D}_2$ and $\mathbf{D}_4$ according to numerical versions of the required conditions. Note that to ensure that the energy remains preserved in the scheme, Eq. (50) must hold, hence both $\mathbf{D}_2$ and $\mathbf{D}_4$ must remain symmetric. Table 1 lists numerical boundary conditions for clamped, simply supported, and free end conditions that satisfy this criterion, making use of the spatial shift operators defined in Section 3.1.

Realistic simulations suitable for sound synthesis require the inclusion of damping terms. Here we restrict ourselves to the losses associated with the vibrating object and its interaction with an external fluid, thus neglecting any impact friction. Such damping can be included in the Euler-Bernoulli model by introducing resistive and Kelvin-Voigt terms [51], rewriting the equation of motion as

$$\rho A \frac{\partial^2 y}{\partial t^2} = \tau \left( \frac{\partial^2 y}{\partial x^2} + \eta \frac{\partial^3 y}{\partial t \partial x^2} \right) - EI \left( \frac{\partial^4 y}{\partial x^4} + \eta \frac{\partial^5 y}{\partial t \partial x^4} \right) - \rho A \gamma \frac{\partial y}{\partial t} + k_c \lfloor (y_c - y)^\alpha \rfloor \quad (56)$$

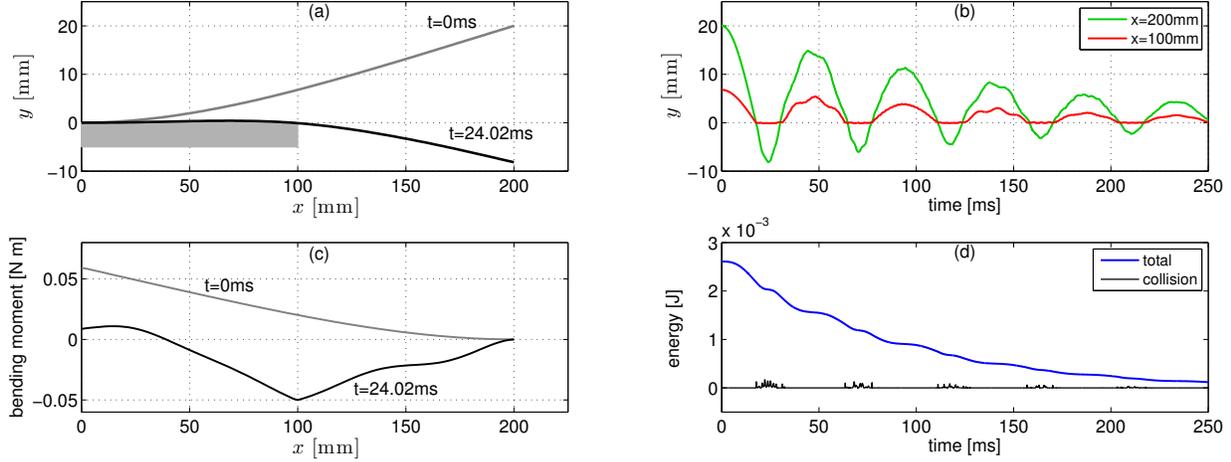

Figure 9: Simulation of a cantilever beam colliding with a flat barrier, for $L = 0.2$ m, $EI = 0.03375$ N m$^2$, $\rho A = 0.03$ kg m$^{-1}$, $\gamma = 10$ s$^{-1}$, $\eta = 10^{-6}$ s, $k_c = 5 \times 10^6$ N m$^{-2}$, $\alpha = 1$, $\Delta x = 1/460$ m, $\Delta t = 1/176.4$ ms. (a) Beam state at two moments in time. (b) Observed displacement at two positions along the beam axis. (c) Bending moment ($M = -EI\partial^2 y/\partial x^2$) at two moments in time. (d) Evolution of energies.

where $\eta$ and $\gamma$ are damping factors which can be loosely associated with internal friction and fluid damping, respectively [10]. Considering now that Eq. (37a) is a direct approximation of Eq. (29a), which is of the dimension of force density, the losses can be included in the numerical formulation simply by adding finite difference approximations of the respective damping terms to Eq. (37a). This procedure leaves the scheme intact except for the nonlinear function, which now takes the form

$$\mathbf{F} = \left[\left(1 + \frac{\gamma \Delta t}{2}\right)\mathbf{I} + \left(1 + \frac{2\eta}{\Delta t}\right)\mathbf{D}\right]\mathbf{s} + 2(\mathbf{D}\mathbf{y}^n - \mathbf{q}^n) + \zeta \mathbf{S}^{-1}\left(\lfloor(\mathbf{y_c} - \mathbf{y}^n - \mathbf{s})^{\alpha+1}\rfloor - \lfloor(\mathbf{y_c} - \mathbf{y}^n)^{\alpha+1}\rfloor\right). \qquad (57)$$

To exemplify these extensions, a cantilever beam colliding with a flat barrier is simulated. The parameters are chosen as listed in the caption of Fig. 9. In this configuration, the model resembles a plastic ruler beating against a flat table. Fig. 9(a) shows the initial state of the beam as well as that after 24.02ms, at which point the tip of the beam reaches its lowest position. Fig. 9(b) demonstrates that the beam is effectively constrained at its halfway point. The correctness of the implementation of the clamped and free boundary conditions can be verified by inspecting Fig. 9(a) and 9(c), the latter displaying the bending moment at the corresponding two instances. The knick point at $t = 24.02$ms indicates the influence of the interaction with the table, which at that moment in time is localised at the table corner point. Fig. 9(d) shows the evolution of the total and the collision energy, the latter defined as $V_c^n = \sum_m \mathcal{V}_c(y_m^n)\Delta x$. This plot reveals a complex bouncing pattern with multiple impacts per cycle, similar to Fig. 8(c). Finally, Fig. 9(b) and 9(d) confirm a damped system behaviour. An energy preservation check does not apply now, but stability may still be observed in that $\partial H/\partial t \leq 0$ at all times.

## 4. Application to the tanpura

The tanpura is a fretless string instrument providing lively sounding drones typical of the musical cultures of the Indian subcontinent. Like various other Eastern string instruments, its specific overtone-rich sound results from the interaction of its strings with a slightly curved bridge, but with the additional feature of having a thin thread placed between the string and the bridge (see Fig. 10(a)), which effectively creates a 'two-point bridge' [52]. Making the simplification of considering the thread and the bridge as immovable objects, the vibrational behaviour can be modelled as a string with simply supported ends meeting a curved bridge placed at a small distance $x_b$ from the thread (see Fig. 10(b)). A more complex model of the tanpura, including impact damping and taking into account the full length of the string (up to the tuning bead) is considered in [53].

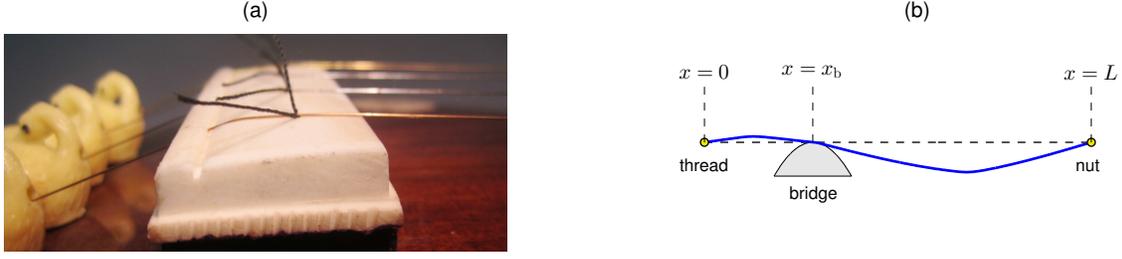

Figure 10: (a) Close-up of a tanpura bridge. (b) Simplified model (dimensions are not proportional).

*4.1. Numerical formulation*

Discretisation of the above model largely proceeds as discussed in Section 3, apart from one aspect. That is, in order to allow specifying a fine spatial detail in the bridge curvature without having to increase the spatial resolution $\Delta x$ of the discretised string, spatial interpolation is applied in the calculation of the contact forces. Setting the spatial resolution of the bridge profile vector $\mathbf{y}_b$ to $\Delta x_b$, an interpolating matrix $\mathcal{I}_b$ is used to translate the string displacements from one grid to the other:

$$\bar{\mathbf{y}}^n = \mathcal{I}_b \mathbf{y}^n, \tag{58}$$

where $\bar{\mathbf{y}}^n$ is a vector holding the interpolated string displacement values; a third-order Lagrange interpolant [43] is applied here. The (scaled) contact forces are formulated at the points on the finer scale, indexed by $i$, as

$$\bar{f}_i^n = \beta_b \frac{\lfloor (y_{b,i} - \bar{y}_i^n - \bar{s}_i)^2 \rfloor - \lfloor (y_{b,i} - \bar{y}_i^n)^2 \rfloor}{\bar{s}_i}, \tag{59}$$

where

$$\beta_b = \frac{k_b \Delta t^2}{2 \rho A}, \quad \bar{s}_i = \bar{y}_i^{n+1} - \bar{y}_i^n \tag{60}$$

The forces can be translated back to the string spatial coordinates using a corresponding downsampling interpolant:

$$\mathbf{f}^n = \mathcal{I}_b^* \bar{\mathbf{f}}^n \tag{61}$$

where, following [8], the downsampling interpolant is defined as the scaled conjugate $\mathcal{I}_b^* = (\Delta x_b/\Delta x)\,\mathcal{I}_b^t$, which ensures that energy conservation is not affected. The numerical formulation remains exactly as before, apart from the non-linear equation to be solved, which now takes the form

$$\mathbf{F} = \left[\left(1 + \frac{\gamma \Delta t}{2}\right)\mathbf{I} + \left(1 + \frac{2\eta}{\Delta t}\right)\mathbf{D}\right]\mathbf{s} + 2\left(\mathbf{D}\mathbf{y}^n - \mathbf{q}^n\right) + \mathbf{f}^n \tag{62}$$

where the force term $\mathbf{f}^n$ is a non-linear function of the 'step vector' $\mathbf{s}$, and where the interpolated version of $\mathbf{s}$ is computed in the same way as for $\mathbf{y}$, i.e. $\bar{\mathbf{s}} = \mathcal{I}_b \mathbf{s}$. This equation can again be solved using the multi-dimensional Newton method, with the Jacobian taking the form

$$\mathbf{J} = \left[\left(1 + \frac{\gamma \Delta t}{2}\right)\mathbf{I} + \left(1 + \frac{2\eta}{\Delta t}\right)\mathbf{D}\right] + \mathcal{I}_b^* \mathbf{G} \mathcal{I}_b, \tag{63}$$

where $\mathbf{G}$ is a diagonal matrix with elements $\{g_{i,i}\} = \partial \bar{f}_i/\partial \bar{s}_i$.

*4.2. Results*

To obtain suitable string parameters, the diameter ($d = 0.3$ mm) and speaking length ($L = 628$ mm) of the third string of a small travelling tanpura were measured. Taking into account the fundamental frequency as well as the mass density and Young's modulus of steel, the tension and stiffness terms were set accordingly to $\tau = 31.47$ N m$^{-1}$ and $EI = 8.35 \times 10^{-5}$ N m$^2$, with $\rho A = 5.58 \times 10^{-4}$ Kg m$^{-1}$. The damping parameters were set to $\gamma = 0.1$ s$^{-1}$ and $\eta = 5 \times 10^{-8}$ s, which results in a frequency-dependent decay pattern that approximately matches that observed when the tanpura string is left in free vibration (i.e. without string-bridge interaction). The bridge contact elasticity coefficient is chosen as $k_b = 5 \times 10^8$ N m$^2$, which ensures that the effective compression does not exceed 1% of the string diameter. The bridge shape is defined as the parabole $\mathbf{y}_b(x) = -4(x_b - x)^2$, which yields a curvature similar to that of the tanpura bridge around the point $x_b = 5$ mm where the string meets the bridge when at rest. The numerical parameters are as follows: $\Delta t = 1/176.4$ ms, $\Delta x = 3.1$ mm, and $\Delta x_b = 0.2$ mm.

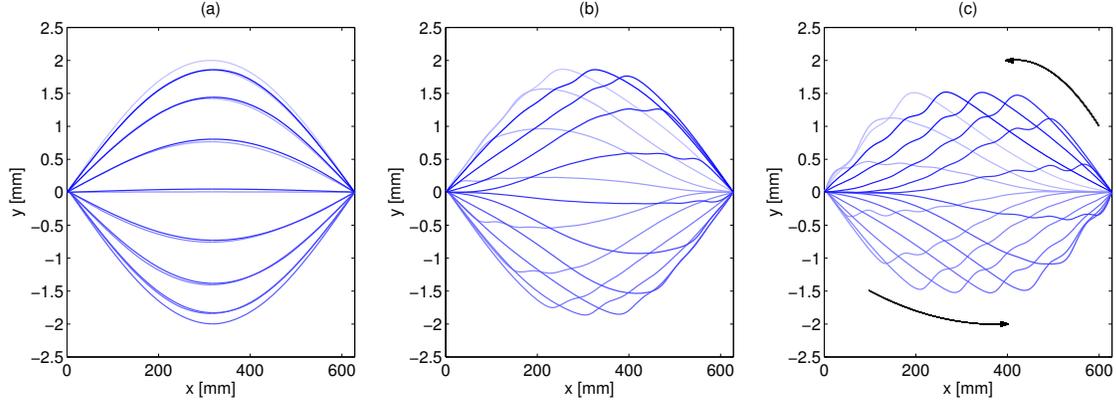

Figure 11: Snapshots of the string motion during the 1st (a), 17th (b) and 33rd (c) period of oscillation. The arrows in (c) indicate the movement of the kink, indicative of a Helmholtz-like motion.

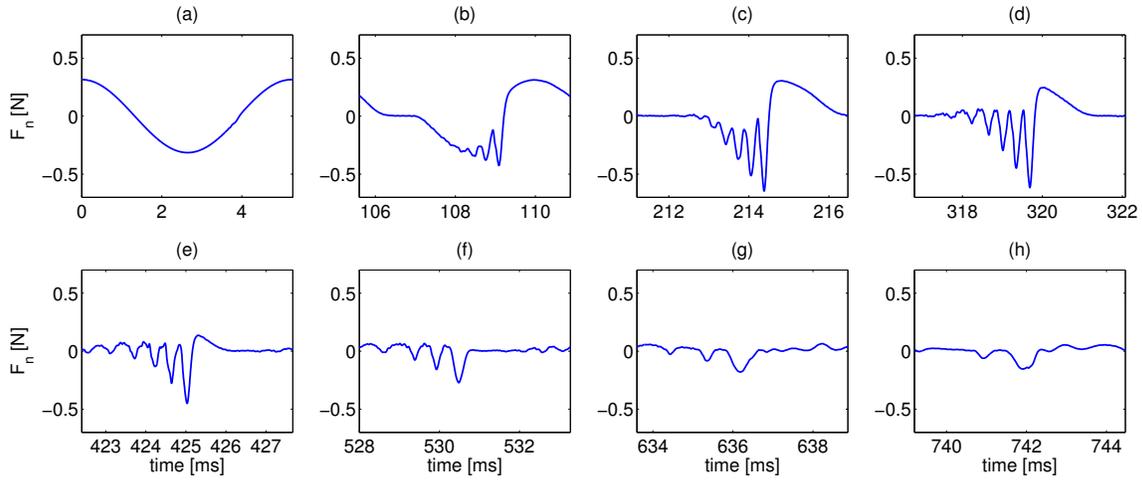

Figure 12: Evolution of the nut force signal. Each plot shows a single period of oscillation. The initial condition was set as $y(x) = 0.002 \sin(\pi x/L)$.

Fig. 11 shows snapshots of the string motion during the 1st, 17th and 33rd period of oscillation, for an initial condition that matches the shape of the first mode of the string. The more recent states are represented by colour-intensive curves, while the colour-tone is fading out for the earlier string states. It can be observed that the bridge interaction forces the string to gradually take on a more triangular shape, indicating the excitation of the other modes of vibration. As such, a Helmholz-like motion emerges, with the kink travelling along the string as indicated by the arrows in Fig. 11(c). Such a motion has been suggested by earlier simulations of the tanpura [54] and has been also encountered in studies of various other string-bridge configurations without the thread [23, 55].

Instances of the corresponding transversal force signal at the nut are shown in Fig. 12, which reveals the gradual development of a precursor in the force waveform. As explained in [52], the precursor is a packet of high-frequencies arriving back at the bridge before the lower frequencies due to the string stiffness. In each string cycle, the precursor

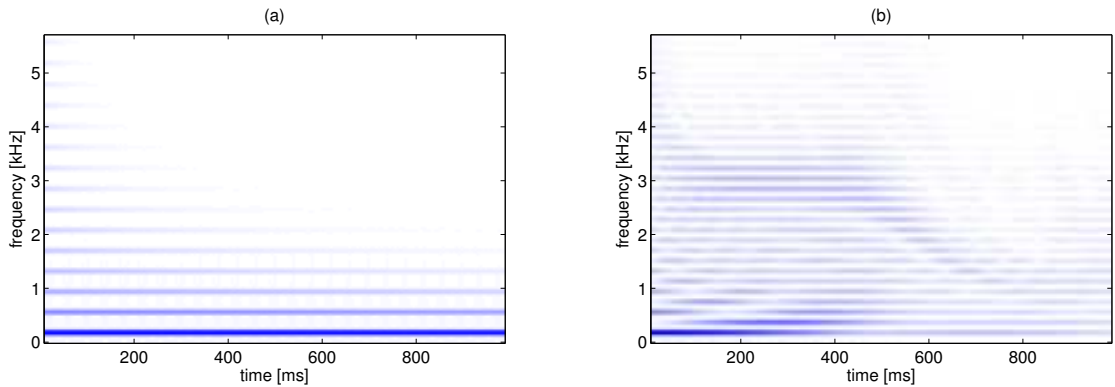

Figure 13: Spectrogram of the nut force signal when plucking the sting at $x = L/2$, for a simulation without bridge interaction (a), and with bridge interaction (b).

is 'fed' with high-frequency components by the nonlinear interaction between the string and bridge.

Fig. 13 shows the spectral evolution of the nut force signal obtained from the simulation when setting the initial condition as a triangular shape that mimics plucking the string at mid-point, with a maximum displacement of 2 mm. Fig. 13(a) confirms that none of the even-numbered string modes are excited in the absence of bridge interaction, and the expected frequency-dependent decay is observed. With the bridge in place (Fig. 13(b)), all modes are excited, and the precursor can be observed as a formant region with a spectral centroid that varies over time. That is, at first the formant frequency decays, then from around $t = 120$ ms to 450 ms it stays approximately constant, followed by a period of slower decay. The appearance of these distinct regimes in the formant frequency decrease pattern are in accordance with the analysis of experimentally obtained nut force signals by Valette et al. [52].

## 5. Conclusion

Numerical schemes for frictionless vibro-impact problems can be derived in stable implicit form by discretising the Hamiltonian description of the system dynamics. Stability follows directly from the inherent energy conservation demonstrated in Sections 2.3 and 3.2. For both lumped and distributed contact the presented formulation involves solving, at each time step, a nonlinear equation which is guaranteed to have a unique solution. The generally excellent convergence conditions for Newton iteration are underlined by the convexity proof for lumped contact given in Section 2.2.

Because of these properties the proposed methodology can be used to formulate time-domain models with improved robustness in comparison to methods previously applied to musical instrument simulation. That is, unlike schemes derived by discretising a Newtonian description, the simulations do not suffer from energy jumps during the decoupling of the impacting objects. As a result, implementations require neither the energy corrections employed in various other formulations [17, 19] nor the monitoring of energy flows common to wave digital filter models [26].

The approach is distinct from that taken in [46] in that it is fully Hamiltonian, and in the absence of damping terms the preserved quantity exactly equals the energy of the underlying continuous model, which is one known way to measure the success of a numerical simulation [34]. Furthermore, there are no stability bounds to respect (i.e. the stability is unconditional). A comparatively positive feature of the approach in [46] is that numerical dispersion can be controlled better, which is advantageous for cases in which aliasing and errors in the effective repelling force are not likely to be prevalent; this advantage is less likely to hold for modelling repetitive contact with nearly rigid barriers, given that numerical dispersion becomes negligible for the sampling frequencies required to avoid aliasing and the accumulation of repelling force errors.

As shown in Section 3.6, for distributed contact the boundary conditions can be altered without affecting energy conservation, and the resonator can be made dissipative simply by adding loss terms. Considering that in Hunt-

Crossley form, the power law for distributed contact can be written as

$$f = k_c \lfloor (y_c - y)^\alpha \rfloor \left(1 - r_c \frac{\partial y}{\partial t}\right) = -\mathcal{V}'_c(y) + r_c \frac{\partial \mathcal{V}_c(y)}{\partial t} \tag{64}$$

the extension to nonlinear impact damping is similarly straightforward.

These results pave the way for improved time-domain modelling of musical instruments for which impactive interaction represents an essential ingredient of the sound production process, especially if the contact is repetitive or distributed. An illustrative example has been presented here, simulating the drone sounds emitted by the tanpura. The generation of high frequency components due to the interaction of the string with a flat, two-point bridge can be simulated using the proposed approach, avoiding numerical artefacts that have appeared in previous studies (see e.g. [54, Fig. 34]). In comparison to the digital waveguide model used in [23, 24], the proposed approach has the advantage that the elasticity properties of the contact—which for some of such instruments are known to play an important role [52]—can be specified. Other impactive interactions of musical interest that the proposed schemes can be directly applied to include reed beating in woodwind instruments [56], string-bridge coupling in pianos [24], and braypin-string collisions in early harps.

The main limitation of the methodology is that the convergence of the Newton iterations, and therefore the invariance of the numerical energy, is subject to finite-precision arithmetic. Namely a bound exists on the contact stiffness parameter $k_c$; extreme values of it can result in there not being sufficient number precision available in updating the displacement step **s** such that the resulting energy step $H^{n+1} - H^n$ falls within machine precision. For example, at standard audio rates the simulation of a string interacting with a barrier characterised by $\alpha = 1$ would typically allow $k_c$ to be set up to a value between $10^7$ and $10^8$ with full Newton convergence. When higher values are needed, a smaller time step may have to be used, coming at the cost of an increase in computational complexity. As discussed in Sections 2.6 and 3.5, some oversampling may be needed anyway in sound synthesis applications in order to avoid the generation of perceptually salient aliasing components.

Given that all mechanical systems can be formulated in Hamiltonian form, the expectation is that the proposed approach naturally extends to other musical instrument phenomena, such as hammer-string and mallet-membrane interaction, as well as to vocal-fold collision in speech and singing. Further noteworthy directions to explore include an alternative formulation of the distributed object, such as a Timoshenko model [57], which would facilitate the inclusion of geometric nonlinearities. Such an extension has in fact recently been shown to be compatible with closely related methods for time-domain modelling of vibrating piano strings [39].

As a final point of interest, it is worth mentioning that the unconditional stability of the proposed schemes is a useful property regarding on-line variation of model parameters, which is desirable from the perspective of studying and modelling player control.